\documentclass{article}
\usepackage{amsmath}
\usepackage{amssymb}
\usepackage{diagrams}
 \newtheorem{thm}{Theorem}[section]
 \newtheorem{cor}[thm]{Corollary}
 
 \newtheorem{claim}[thm]{Claim}
 \newtheorem{prop}[thm]{Proposition}
 \newtheorem{exer}[thm]{Exercise}
 \newtheorem{conv}[thm]{Convention}

 \newtheorem{rem}[thm]{Remark}
 \newtheorem{nota}[thm]{Notation}

 \newtheorem{defn}[thm]{Definition}
 \newtheorem{ex}[thm]{Example}
 \numberwithin{equation}{section}
\newcommand{\spec}{{\rm Spec}}
\newcommand{\gr}{{\rm grade}}
\newcommand{\de}{{\rm depth}}

\newcommand{\supp}{{\rm Supp}}
\newcommand{\hgt}{{\rm ht}}

\newcommand{\fd}{{\rm fd}}
\newcommand{\reg}{{\rm Reg}}
\newcommand{\cmd}{{\rm cmd}}
\newcommand{\smin}{{\rm Min}}
\newcommand{\ass}{{\rm Ass}}

\begin{document}  
  \title{About the Cohen-Macaulay defect and almost Cohen-Macaulay rings}
  \author{Cristodor Ionescu }
  \date{}
\maketitle

\begin{abstract}
We notice the connection between almost Cohen-Macaulay rings and the Cohen-Macaulay defect. We introduce a Serre-type condition for modules that is connected to the Cohen-Macaulay defect in the same way that the condition $(S_n)$ is connected to Cohen-Macaulay modules.
\end{abstract}
\maketitle
\vspace{5mm}
\par All rings considered will be commutative with unit and Noetherian. All modules are  finitely generated.

\section{Introduction and motivation}

\par Almost Cohen-Macaulay rings have shown up due to a flaw in the first edition of Matsumura's book \cite{Mat}:
\begin{claim}\label{matflaw} \cite[p.97, 15.C]{Mat} If A is an
arbitrary Noetherian ring and $\mathfrak{p}\in \spec(A),$ we have
$\de(M_\mathfrak{p}) = 0$ as $A_\mathfrak{p}$-module $\Leftrightarrow \mathfrak{p}A_\mathfrak{p}\in\ass_{A_\mathfrak{p}}(M_\mathfrak{p})\Leftrightarrow \mathfrak{p}\in\ass_A(M)$.
It follows that, in general, $\de(M_\mathfrak{p})$ as $A_\mathfrak{p}$-module is equal
to $\de_\mathfrak{p}(M)$.
\end{claim}
\par\noindent The error was corrected in the second edition of the book \cite{Mat2}:
\begin{exer}\label{matcor} \cite[p. 113, Ex. 1]{Mat2}  Find an example of a Noetherian local ring $A$ and a finite $A$-module $M$ such that $\de(M)>\de(A).$ Also, find $A, M$ and $\mathfrak{p}\in\spec(A)$ such that $\de(M_\mathfrak{p})>\de_\mathfrak{p}(M)$.
\end{exer}

\par More concretely in \cite{Mat3} the following exercise is proposed to the reader:
\begin{exer}\label{mat3} \cite[Ex. 16.5]{Mat3}
Let $A$ be a Noetherian local ring, $M$ a finite $A$-module, and $\mathfrak{p}$ a prime ideal of $A$; show that $\de_\mathfrak{p}(M)\leq\de_{A_\mathfrak{p}}(M_\mathfrak{p})$ and construct an example where the inequality is strict.
\end{exer}
In the section of \cite{Mat3} dedicated to \textit{Solutions and hints for the exercises}, the solution to Exercise \ref{mat3}  shows that if $k$ is a field and $A=(k[X,Y,Z]/(X,Y,Z)^2\cap(Z))_{(X,Y,Z)}$,  then $\de(A)=0$ and $\de(A_\mathfrak{p})>0,$ where   $\mathfrak{p}=(X,Z)A.$
\par In 1998, in a paper written in Chinese  \cite{Han},   Han studied the class of commutative Noetherian rings $A$ such that for any prime ideal $\mathfrak{p}$ in $A,$ we have $\de_\mathfrak{p}(A)=\de_{\mathfrak{p}A_\mathfrak{p}}(A_\mathfrak{p})$ and called them D-rings.
Afterwards, Kang \cite[Def. 1.2]{K1}  renamed these rings almost Cohen-Macaulay rings, proved several properties of almost Cohen-Macaulay rings and introduced also almost Cohen-Macaulay modules.
\begin{defn}\label{acmmod1} {\rm \cite[Def. 1.2]{K1}}
 A finitely generated $A$-module M is called almost Cohen-Macaulay if $\de_\mathfrak{p}(M)=\de_{\mathfrak{p}A_\mathfrak{p}}(M_\mathfrak{p})$ for any prime ideal $\mathfrak{p}\in\supp(M).$
\end{defn}
After the paper \cite{K1} was published, several papers  about almost Cohen-Macaulay modules and  rings were written by several authors: \cite{CTT}, \cite{I},  \cite{IT} \cite{K2},  \cite{MT}, \cite{MT2} etc.
\par It follows from \cite[Lemma 2.6]{K1} that a Noetherian local ring $A$ is almost Cohen-Macaulay if and only if $\dim(A)-\de(A)\leq 1.$ A similar statement is valid for modules, as follows from \cite[Lemma 2.6]{K1}.
\par It turns out that if $M$ is a finitely generated $A$-module, the integer $\dim(M)-\de(M)$ occuring in the above relation, was already considered a long time ago. It was mentioned by Grothendieck \cite[D\'ef. 16.4.9]{EGA20} and called the codepth (coprofondeur) of $M.$ It was renamed the {\it Cohen-Macaulay defect} of $M$ and denoted $\cmd(M)$ in more recent papers (see e.g. \cite{AF1} and \cite{AFH}). It also turns out that the rings $A$ with the property that $\cmd(A)\leq 1$ and even more generally $\cmd(A)\leq n$ for some natural number $n,$ were roughly considered by Grothendieck \cite[\S\S \  5,6,7]{EGA}. It is the first aim of the present paper to give a  unitary view and treatment of these objects and to recall some of the results spread in the literature. The second aim is to introduce a new condition, called $(C_n^l),$ condition that plays a role similar to the well-known property of Serre. The main properties of this condition  and the connections between this property and the Cohen-Macaulay defect will be investigated, first in the case of rings and then in the case of finitely generated $A$-modules.

\section{Cohen Macaulay defect and almost Cohen-Macaulay rings}
In this section we will mention, with precise citations but with no proof, the main notions and results concerning the Cohen-Macaulay defect that are showing the departure point of the present paper.
\begin{defn}\label{cmdef}{\rm \cite[D\'ef.\!\! 16.4.9]{EGA20}, \cite[D\'ef.\! 5.7.12]{EGA}, \cite[Introduction]{AF1}} Let $A$ be a  Noetherian local ring and $M$ a finitely generated $A$-module. Then the non-negative integer $\cmd_A(M):=\dim(M)-\de_A(M)$ is called the Cohen-Macaulay defect of $M$. If $M=(0)$, we put $\cmd_A(M)=0.$ If $A$ is not local, then the Cohen-Macaulay defect of $M$ is $$\cmd_A(M)=\sup\{\cmd_{A_\mathfrak{p}}(M_\mathfrak{p})\ \vert\ \mathfrak{p}\in\spec(A)\}.$$
\end{defn}

\begin{defn}\label{acmdef} A Noetherian  ring $A$ is called an almost Cohen-Macaulay ring if $\cmd_A(A)\leq 1.$
\end{defn}
\begin{defn}\label{acmmod} If $A$ is a Noetherian ring, a finitely generated $A$-module $M$ is called almost Cohen-Macaulay if $\cmd_A(M)\leq 1.$
\end{defn}
\begin{rem}\label{dim1} i) Any Cohen-Macaulay module is almost Cohen-Macaulay.
\par ii) Any Noetherian local ring $A$ with $ dim(A)\leq 1$ is almost Cohen-Macaulay.

\end{rem}
\begin{ex}\label{excm}
Let $k$ be  a field and $0\leq d<n$ be natural numbers. Let $r=n-d$. If $d=0$ let $S:=k[X_0,\ldots,X_r]$ and if $d>0$ let $S:=k[X_0,\ldots,X_r,T_1,\ldots,T_d].$ Let $I:=(X_0)\cap(X_0,\ldots,X_r)^{r+1}$ and $A:=(S/I)_{(X_0,\ldots,X_r,T_0,\ldots,T_d).}$ Then $\dim(A)=r+d=n$ and $\de(A)=d$, that is $\cmd(A)=r.$ This shows that there are local rings of any Cohen-Macaulay defect.
\end{ex}
Below is a list of known facts about Cohen-Macaulay defect and almost Cohen-Macaulay rings.

\begin{prop}\label{regel} {\rm \cite[Prop.\! 16.4.10]{EGA20}, \cite[Lemmas 2.6, 2.7]{K1}}
 Let $(A,\mathfrak{m})$ be a Noetherian local ring, M a finitely generated $A$-module and $x\in\mathfrak{m}$ a regular element. Then $\cmd_A(M)=\cmd_{A/xA}(M/xM).$ In particular, M is an almost Cohen-Macaulay  $A$-module if and only if $M/xM$ is an almost Cohen-Macaulay  $A/xA$-module.
\end{prop}
\begin{prop}\label{compl}{\rm \cite[Prop. 16.4.10]{EGA20}, \cite[Cor. 2.4]{I}} Let $A$ be a Noetherian local ring and $M$ a finitely generated $A$-module. Then $\cmd_A(M)=\cmd_{\widehat{A}}(\widehat{M}).$ In particular, $M$ is an almost Cohen-Macaulay $A$-module  if and only if $\widehat{M}$ is an almost Cohen-Macaulay $\widehat{A}$-module.
\end{prop}

\begin{prop}\label{local} {\rm \cite[Prop. 6.11.5]{EGA}, \cite[Lemma 2.6]{K1}}
 Let $(A,\mathfrak{m})$ be a Noetherian local ring and M a finitely generated $A$-module. Then $\cmd_{A_\mathfrak{p}}(M_\mathfrak{p})\leq\cmd_A(M)$ for any prime ideal $\mathfrak{p}\in\spec(A).$ In particular, $M$ is an almost Cohen-Macaulay $A$-module if and only if  $M_\mathfrak{p}$ is an almost Cohen-Macaulay $A_\mathfrak{p}$-module for any $\mathfrak{p}\in\spec(A).$
\end{prop}

\begin{prop}\label{flcmd} {\rm \cite[Cor. 2.5]{BCLM}}
Let $A$ be a commutative ring, and let $B$ and $C$ be $A$-algebras such that $B\otimes_AC$ is Noetherian. Let $\mathfrak{r}\in\spec(B\otimes_AC)$ and set $\mathfrak{p}=\mathfrak{r}\cap B, \mathfrak{q}=\mathfrak{r}\cap C$ and $\mathfrak{n}=\mathfrak{r}\cap A,$ so that we have the following commutative diagram:
\begin{diagram}[size=2.1em]
A_\mathfrak{n}   &     \rTo   &      B_\mathfrak{p} \\
\dTo &           & \dTo  \\
C_\mathfrak{q}    &\rTo      & (B\otimes_AC)_\mathfrak{r}
\end{diagram}
Assume that $B_\mathfrak{p}$ is flat over $A_\mathfrak{n}.$
Then
$$\cmd((B\otimes_AC)_\mathfrak{r})=\cmd(B_\mathfrak{p})+\cmd(C_\mathfrak{q})-\cmd(A_\mathfrak{n})=$$
$$=\cmd(B_\mathfrak{p}/\mathfrak{n}B_\mathfrak{p})+\cmd(C_\mathfrak{q}).$$
\end{prop}

\begin{cor}\label{consec}
 Let $A$ be a commutative ring and let $B$ and $C$ be $A$-algebras such that $B\otimes_AC$ is Noetherian. Let $\mathfrak{r}\in\spec(B\otimes_AC)$ and set $\mathfrak{p}=\mathfrak{r}\cap B, \mathfrak{q}=\mathfrak{r}\cap C$ and $\mathfrak{n}=\mathfrak{r}\cap A.$ Assume that $B_\mathfrak{p}$ is flat over $A_\mathfrak{n}$ and that $A_\mathfrak{n}$ is Cohen-Macaulay. If $B_\mathfrak{p}$ and $C_\mathfrak{q}$ are almost Cohen-Macaulay and one of them is Cohen-Macaulay, then $(B\otimes_AC)_\mathfrak{r}$ is almost Cohen-Macaulay.
\end{cor}

\begin{cor}\label{fieldext}{\rm \cite[Prop. 6.7.1]{EGA}, \cite[Cor. 2.12]{I}} Let $k$ be a field, A a Noetherian $k$-algebra and $L$ a field extension of $k$ such that $B:=A\otimes_kL$ is Noetherian. Let $\mathfrak{q}\in\spec(B)$ and $\mathfrak{p}=\mathfrak{q}\cap A.$ Then $\cmd(B_\mathfrak{q})\leq\cmd(A_\mathfrak{p}).$ In particular, if $A$ is almost Cohen-Macaulay, then $A\otimes_kL$ is almost Cohen-Macaulay.
\end{cor}

\begin{prop}\label{morph}{\rm \cite[Prop.\! 6.3.2]{EGA}, \cite[Th.\! 2.10]{K1}, \cite[Prop. 2.2]{I}}
Let $u:(A,\mathfrak{m},k)\to B$ be a morphism of Noetherian local rings, $M\neq(0)$ a finitely generated $A$-module and $N\neq(0)$ a finitely generated $B$-module which is $A$-flat. Then
$$\cmd_B(M\otimes_AN)=\cmd_A(M)+\cmd_{B\otimes_Ak}(N\otimes_Ak).$$
In particular, if $u$ is flat then
$$\cmd(B)=\cmd(A)+\cmd(B/\mathfrak{m}B).$$
\end{prop}

\begin{cor}\label{flat}{\rm \cite[Prop. 2.2]{I}}
 Let $u:(A,\mathfrak{m},k)\to B$ be a flat morphism of Noetherian local rings. Then:
 \par\noindent i) If $B$ is almost Cohen-Macaulay, then $A$ and $B/\mathfrak{m}B$ are almost Cohen-Macaulay.
 \par\noindent ii) If $A$ and $B/\mathfrak{m}B$ are almost Cohen-Macaulay and one of them is Cohen-Macaulay, then $B$ is almost Cohen-Macaulay.
\end{cor}
More generally:
\begin{prop}\label{fflat}{\rm\cite[Th. 3.1]{AFH}}
Let $u:A\to B$ be a morphism of Noetherian local rings such that $\fd_A(B)<\infty.$ Then $\cmd(A)\leq\cmd(B).$
\end{prop}

\begin{prop}\label{open}{\rm \cite[Prop. 6.11.2]{EGA}, \cite[Cor. 3.10]{I}} Let $A$ be a Noetherian local ring that is a quotient of a regular local ring and $n$ a natural number. Then the set
$$\cmd_n(A)=\{\mathfrak{p}\in\spec(A)\ \vert \ \cmd(A_\mathfrak{p})\leq n\}$$ is open in $\spec(A).$ In particular, the set
$$aCM(A):=\{\mathfrak{p}\in\spec(A)\ \vert \ A_\mathfrak{p}\ \text{is almost Cohen-Macaulay}\}=\cmd_1(A)$$ is open in $\spec(A).$
\end{prop}

\section{Serre-type  properties for rings}

Recall the following well-know definition, usually attributed  to Serre:
\begin{defn}\label{sn} {\rm \cite[D\'ef. 5.7.2]{EGA} }
Given a natural number $n,$ a Noetherian ring $A$ is said to have Serre's property $(S_n)$ if $\de(A_\mathfrak{p})\geq\min(\hgt({\mathfrak{p}}),n)$ for any prime ideal $\mathfrak{p}\in\spec(A).$
\end{defn}
\begin{rem}\label{propsn}
 It is easy to see that $(S_n)\Rightarrow(S_{n-1})$ and that the property $(S_n)$ localizes well, in the sense that if $A$ has the property $(S_n)$ and $\mathfrak{p}$ is a prime ideal of $A$, then $A_\mathfrak{p}$ has the property $(S_n).$
\end{rem}
Moreover, the following facts show the connection between the properties $(S_n)$ and Cohen-Macaulay:

\begin{rem}\label{cmsn}{\rm \cite[(17.I)]{Mat}} A Noetherian ring $A$ is Cohen-Macaulay if and only if $A$ has the property $(S_n)$ for any $n\in\mathbb{N}$.
\end{rem}

\begin{rem}\label{sncm}
  A Noetherian ring $A$ has the property $(S_n)$ if and only if $A_\mathfrak{p}$ is Cohen-Macaulay for any $\mathfrak{p}\in\spec(A)$ with $\de(A_\mathfrak{p})<n.$
\end{rem}

\par\noindent\textit{Proof:} Assume that $A$ has $(S_n)$ and let $\mathfrak{p}\in\spec(A)$ with $\de(A_\mathfrak{p})<n.$ By $(S_n)$ we have that $\de(A_\mathfrak{p})\geq\min(n,\hgt(\mathfrak{p})).$ From this it follows that $\hgt(\mathfrak{p})<n$ and $\de(A_\mathfrak{p})\geq\hgt(\mathfrak{p}),$ consequently $A_\mathfrak{p}$ is Cohen-Macaulay.
\par\noindent Conversely, assume that $A_\mathfrak{p}$ is Cohen-Macaulay for any $\mathfrak{p}\in\spec(A)$ with $\de(A_\mathfrak{p})<n.$ Let $\mathfrak{p}\in\spec(A)$ with $\de(A_\mathfrak{p})<n.$ Then $n>\hgt(\mathfrak{p})=\de(A_\mathfrak{p})=\min(n,\hgt(\mathfrak{p})).$ If $\de(A_\mathfrak{p})\geq n,$ it is clear that $\de(A_\mathfrak{p})\geq\min(n,\hgt(\mathfrak{p})$.

 In \cite{I} the property $(S_n)$ was generalized, in order to be useful in the characterization of almost Cohen-Macaulay rings.

 \begin{defn}\label{acmcn} {\rm \cite[Def. 3.1]{I}} Given a natural number $n,$ a Noetherian ring $A$ is said to have  the property $(C_n)$  if $\de(A_\mathfrak{p})\geq\min(n,\hgt(\mathfrak{p}))-1,$ for any $\mathfrak{p}\in\spec(A).$
 \end{defn}

 Then, in \cite[Th. 3.3, Prop. 3.4]{I}, statements similar to Remarks \ref{cmsn} and \ref{sncm} above were proved.
We try to go a little bit further in this direction, extending this idea to rings with Cohen-Macaulay defect less than a given number. Thus we introduce the following notion.

\begin{defn}\label{tn}
Let $n$ and $l$ be natural numbers. We say that a Noetherian ring $A$ has the property $(C_n^l)$ if $$\de(A_\mathfrak{p})\geq\min(\hgt(\mathfrak{p}),n)-l, \forall\ \mathfrak{p}\in\spec(A).$$
\end{defn}
\begin{rem}\label{relations}
\par i) It is obvious that $(C_n^1)$ coincides with $(C_n)$ from \cite{I}.
\par ii) Actually, if $n\geq \dim(A),$ then $\min(\hgt(\mathfrak{p}),n)=\hgt(\mathfrak{p})$, hence the ring $A$ has the property $(C_n^l)$ if and only if $\cmd(A)\leq l.$
\par ii) One can see at once that
$$A\ \text{has the property}\ (C_n^l)\Leftrightarrow\begin{cases}
  \cmd(A_\mathfrak{p})  &  \text{if}\  \hgt(\mathfrak{p})<n \\
  \de(A_\mathfrak{p})\geq n-l &  \text{if}\ \hgt(\mathfrak{p})\leq n
\end{cases}$$
\par iii) Note that if  $l\geq\min(\hgt(\mathfrak{p}),n)$ for any prime ideal $\mathfrak{p},$  then $(C_n^l)$  obviously holds. This happens for example if $l>\dim(A).$

\par iv) One can see immediately that $(C_n^l)\Rightarrow (C_{n-1}^l)$ and that  $(C_n^l)\Rightarrow (C_{n}^{l+1}).$
\par v) It is also clear that the property $(C_n^l)$ localizes well,  that is if $A$ has $(C_n^l),$ then $A_\mathfrak{p}$ has $(C_n^l)$ for any prime ideal $\mathfrak{p}\in\spec(A).$
\end{rem}
The following result shows that  $(C_n^l)$ is connected to the property $\cmd(A)\leq l$ in the same way as  $(S_n)$ is connected to the Cohen-Macaulay property.


\begin{thm}\label{charcmd} Let A be a Noetherian ring and $n,l$ be natural numbers. Then $\cmd(A)\leq l$ if and only if $A$ has $(C_n^l)$ for any $n\in\mathbb{N}.$
\end{thm}
\par\noindent\textit{Proof:} Assume that $\cmd(A)\leq l$ and let $\mathfrak{p}\in\spec(A).$ Then $\cmd(A_\mathfrak{p})\leq l,$ that is $\hgt(\mathfrak{p})-\de(A_\mathfrak{p})\leq l.$ If $\hgt(\mathfrak{p})\leq n,$ then $\min(n,\hgt(\mathfrak{p}))=\hgt(\mathfrak{p})),$ hence $\min(\hgt(\mathfrak{p},n)-\de(A_\mathfrak{p})=\hgt(\mathfrak{p})-\de(A_\mathfrak{p})\leq l$. If $\hgt(\mathfrak{p})>n,$ then $\min(n,\hgt(\mathfrak{p}))=n,$ hence $\de(A_\mathfrak{p})\geq\hgt(\mathfrak{p})-l\geq n-l=\min(\hgt(\mathfrak{p}),n)-l.$
\par\noindent Conversely, let $\mathfrak{p}\in\spec(A), \hgt(\mathfrak{p})=s.$ Then
$$\de(A_\mathfrak{p})\geq\min(s,\hgt(\mathfrak{p})
)-l=\hgt(\mathfrak{p})-l,$$
that is $\cmd(A_\mathfrak{p})=\hgt(\mathfrak{p})-\de(A_\mathfrak{p})\leq l.$

\begin{prop}\label{charcnk} Let A be a Noetherian ring and $n,l$ be natural numbers. Then  $A$ has $(C_n^l)$ if and only if $\cmd(A_\mathfrak{p})\leq l$ for any prime ideal $\mathfrak{p}$ such that $\de(A_\mathfrak{p})\leq n-l-1.$
\end{prop}

\par\noindent\textit{Proof:} Assume that $A$ has the property $(C_n^l).$ Let $\mathfrak{p}\in\spec(A)$ be such that $\de(A_\mathfrak{p})\leq n-l-1.$ Hence $\min(n,\hgt(\mathfrak{p}))-l\leq\de(A_\mathfrak{p})\leq n-l-1.$
\par\noindent If $\hgt(\mathfrak{p})\leq n,$ then $\hgt(\mathfrak{p})-l\leq\de(A_\mathfrak{p}),$ that is $\cmd(A_\mathfrak{p})\leq l.$
\par\noindent If $\hgt(\mathfrak{p})>n,$ then $n-l\leq\de(A_\mathfrak{p})\leq n-l-1,$ but this is clearly a  contradiction.
\par Conversely, let $\mathfrak{p}\in\spec(A).$
\par\noindent If $\de(A_\mathfrak{p})\leq n-l-1,$ then $\cmd(A_\mathfrak{p})=\hgt(\mathfrak{p})-\de(A_\mathfrak{p})\leq l,$ hence $\hgt(\mathfrak{p})\leq n-1$. This means that $\min(n,\hgt(\mathfrak{p}))=\hgt(\mathfrak{p}),$ that is $\de(A_\mathfrak{p})\leq\min(n,\hgt(\mathfrak{p}))-l.$
\par\noindent If $\de(A_\mathfrak{p})> n-l-1,$ then $\hgt(\mathfrak{p})>n-2,$ hence $\de(A_\mathfrak{p})\geq\min(n,\hgt(\mathfrak{p}))-l.$

\begin{prop}\label{regelcnk}
Let $A$ be a Noetherian ring, $n,l\in\mathbb{N},$ and let $x\in A$ be a non zero divisor. If $A/xA$ has the property $(C_n^l)$, then $A$ has the property $(C^l_n).$
\end{prop}

\par\noindent\textit{Proof:}
Let $\mathfrak{q}\in\spec(A)$ such that $\de(A_\mathfrak{q})=s\leq n-l-1.$ If $x\in \mathfrak{q},$ then $\de(A/xA)_\mathfrak{q}=s-1\leq n-l-2.$ Then $\hgt(\mathfrak{q}/xA)\leq s-1+1=s,$
hence $\hgt(\mathfrak{q})\leq s+1=\de (A_\mathfrak{q})+1.$  If $x\notin \mathfrak{q},$ let $\mathfrak{p}\in\smin(\mathfrak{q}+xA).$ Then $(\mathfrak{q}+xA)A_\mathfrak{p}$ is $\mathfrak{p}A_\mathfrak{p}$-primary and $\de(A_\mathfrak{p})\leq\de(A_\mathfrak{q})+1=s+1.$
Then $\de(A/xA)_\mathfrak{q}=s-1,$ hence $\hgt(\mathfrak{p}/xA)\leq s.$ It follows that $\hgt(\mathfrak{p})\leq s+1=\de(A_\mathfrak{p})+1.$

\begin{nota} In the following we shall consider a property \textbf{P}  of Noetherian local rings, like for example regular, Gorenstein, Cohen-Macaulay etc. For a Noetherian ring $A$ and such a property \textbf{P}, we denote by $\textbf{P}(A)$ the \textbf{P}-locus of $A$, that is $\textbf{P}(A):=\{\mathfrak{p}\in\spec(A)\ \vert\ A_\mathfrak{p}\ \text{has the property \textbf{P}}\}$(cf. \cite[p.187]{Mat3}).
\end{nota}
\par Let us remind the following definition:

\begin{defn}\label{nagcr} {\rm \cite[p. 187]{Mat3}}
We say that a property \textbf{P} of Noetherian local rings satisfies Nagata's Criterion (NC) if the following holds:
if $A$ is a Noetherian ring such for every $\mathfrak{p}\in \textbf{P}(A)$ the set $\textbf{P}(A/\mathfrak{p})$ contains a non-empty open set of $\spec(A/\mathfrak{p})$, then $\textbf{P}(A)$ is open in the Zariski topology of $\spec(A).$
\end{defn}

\begin{thm}\label{Tknagata}
Let $n,l\in\mathbb{N}.$ The property $(C_n^l)$ satisfies (NC).
\end{thm}

\par\noindent\textit{Proof:}
Let $\mathfrak{q}\in C_n^l(A).$ Then $\de(A_\mathfrak{q})\geq\min(n,\hgt(\mathfrak{q}))-l.$
\par\noindent Case a): $\hgt(\mathfrak{q})\leq n.$ Then $\min(n,\hgt(\mathfrak{q}))=\hgt(\mathfrak{q}),$ hence $\de(A_\mathfrak{q})\geq\hgt(\mathfrak{q})-l$, that is $\cmd(A_\mathfrak{q})\leq l.$  By \cite[Prop. 6.10.6]{EGA}, there exists an element $f\in A\setminus \mathfrak{q}$ such that
$$\dim(A_\mathfrak{p})=\dim(A_\mathfrak{q})+\dim(A_\mathfrak{p}/\mathfrak{q}A_\mathfrak{p})$$
 and
 $$\de(A_\mathfrak{p})=\de(A_\mathfrak{q})+\de(A_\mathfrak{p}/\mathfrak{q}A_\mathfrak{p})$$
for any $\mathfrak{p}\in D(f)\cap V(\mathfrak{q})\cap NC_n^l(A).$ Then $\de(A_\mathfrak{p})\ngeqq\min(n,\hgt(\mathfrak{p}))-l.$
\par\noindent Case a1): $\hgt(\mathfrak{p})\leq n.$ Then $\min(n,\hgt(\mathfrak{p}))=\hgt(\mathfrak{p}),$ hence $\de(A_\mathfrak{p})+l<\hgt(\mathfrak{p}).$ Then $$\de(A_\mathfrak{p}/\mathfrak{q}A_\mathfrak{p})+l=\de(A_\mathfrak{p})-\de(A_\mathfrak{q})+l<$$
$$<\hgt(\mathfrak{p})-\de(A_\mathfrak{q})\leq\hgt(\mathfrak{p})-\hgt(\mathfrak{q})+l.$$
Then $\de(A_\mathfrak{p}/\mathfrak{q}A_\mathfrak{p})<\dim(A_\mathfrak{p}/\mathfrak{q}A_\mathfrak{p})=\dim(A_\mathfrak{p})-\dim(A_\mathfrak{q}),$ and it follows that $A_\mathfrak{p}/\mathfrak{q}A_\mathfrak{p}$ is not $(C_n^l).$
\par\noindent Case a2): $\hgt(\mathfrak{p})> n.$ Then $\min(n,\hgt(\mathfrak{p}))=n,$ hence $\de(A_\mathfrak{p})<n-l.$ It follows that
$$\de(A_\mathfrak{p}/\mathfrak{q}A_\mathfrak{p})=\de(A_\mathfrak{p})-\de(A_\mathfrak{q})<n-l-\hgt(\mathfrak{q}).$$
This implies that $A_\mathfrak{p}/\mathfrak{q}A_\mathfrak{p}$ is not $(C_n^l).$
\par\noindent Case b): $\hgt(\mathfrak{q})>n.$ Then $\min(n,\hgt(\mathfrak{q}))=n$ and $\de(A_\mathfrak{q})+l\geq n.$ Since $\hgt(\mathfrak{p})>n,$ it follows that $\min(n,\hgt(\mathfrak{p}))=n$ and $\de(A_\mathfrak{p})+1<n.$ Let $x_1,\ldots,x_r$ be an
$A_\mathfrak{q}$-regular sequence. Then there exists $f\in A\setminus \mathfrak{q}$ such that $x_1,\ldots,x_r$ is $A_f$-regular. If $\mathfrak{p}\in D(f)\cap V(\mathfrak{q}),$ it follows that $A_\mathfrak{p}$ is $(C_n^l).$

\begin{prop}\label{qeopen} Let A be a quasi-excellent ring and $n,l\in\mathbb{N}$. Then $C_n^l (A)$ is open in the Zariski topology of $\spec (A).$
\end{prop}

\par\noindent\textit{Proof:}
Let $\mathfrak{p}\in\spec(A).$ Then $ C_n^l(A/\mathfrak{p})$ contains
the non-empty open set $\reg(A/\mathfrak{p})=\{\mathfrak{p}\in\spec(A)\ \vert\ A_\mathfrak{p}\ \text{is regular}.\}$
Now apply \ref{Tknagata}.

\begin{cor}\label{compopen} Let A be a semilocal complete ring and $l,n\in\mathbb{N}$. Then $C_n^l (A)$ is open in the Zariski topology of $\spec (A).$
\end{cor}

\begin{prop}\label{flattk}
Let $u:A\to B$ be a flat morphism of Noetherian rings and $l,n\in\mathbb{N}.$ If $B$ has $(C_n^l),$ then $A$ has $(C_n^l).$
\end{prop}

\par\noindent\textit{Proof:}
Using \ref{relations},v) it is clear that we may assume that $A$ and $B$ are local rings and that $u$ is a local morphism. Let $\mathfrak{p}\in\spec(A)$ and $\mathfrak{q}\in\smin(\mathfrak{p}B).$ Then $\dim(B_\mathfrak{q}/\mathfrak{p}B_\mathfrak{q})=0,$ hence
$$\de(A_\mathfrak{p})=\de(B_\mathfrak{q})\geq\min(n,\dim(B_\mathfrak{q}))-l=$$
$$=\min(n,\dim(A_\mathfrak{p}))-l.$$

\begin{prop}\label{flattk2}
Let $u:A\to B$ be a flat morphism of Noetherian rings and $l,n\in\mathbb{N}.$
\par\noindent a) If $A$ has $(C_n^l)$ and all the fibers of $u$ have $(S_n),$ then $B$ has $(C_n^l).$
\par\noindent b) If $A$ has $(S_n)$ and all the fibers of $u$ have $(C_n^l),$ then $B$ has $(C_n^l).$
\end{prop}

\par\noindent\textit{Proof:}
a) Let $\mathfrak{q}\in\spec(B)$ and $\mathfrak{p}=\mathfrak{q}\cap A.$ Then by flatness we have
$$\dim(B_\mathfrak{q})=\dim(A_\mathfrak{p})+\dim(B_\mathfrak{q}/\mathfrak{p}B_\mathfrak{q})$$
and
$$\de(B_\mathfrak{q})=\de(A_\mathfrak{p})+\de(B_\mathfrak{q}/\mathfrak{p}B_\mathfrak{q}).$$
By assumption we have
$$\de(A_\mathfrak{p})\geq\min(n,\hgt(\mathfrak{p}))-l,$$
$$\de(B_\mathfrak{q}/\mathfrak{p}B_\mathfrak{q})\geq\min(n,\dim(B_\mathfrak{q}/\mathfrak{p}B_\mathfrak{q}).$$
Hence we have
$$\de(B_\mathfrak{q})=\de(A_\mathfrak{p})+\de(B_\mathfrak{q}/\mathfrak{p}B_\mathfrak{q})\geq$$
$$\geq\min(n,\hgt(\mathfrak{p}))-l+\min(n,\dim(B_\mathfrak{q}/\mathfrak{p}B_\mathfrak{q}))\geq$$
$$\geq\min(n,\hgt(B_\mathfrak{q})+\dim(B_\mathfrak{q}/\mathfrak{p}B_\mathfrak{q}))=\min(n,\hgt(B_\mathfrak{q}))-l.$$
\par\noindent b) The proof is along the same lines.

Another generalization of Serre's condition $(S_n)$ was considered in \cite{Ho}:

\begin{defn}\label{tnholmes} {\rm \cite[Def. 1.2]{Ho}}
Let $n, l\in\mathbb{N}$ be natural numbers. We say that a Noetherian ring $A$ has the property $(S_n^l)$ if $$\de(A_\mathfrak{p})\geq\min(\hgt(\mathfrak{p}),n-l),\ \text{for all}\  \mathfrak{p}\in\spec(A).$$
\end{defn}



\begin{cor}\label{snand snl}
Let A be a Noetherian ring having the property $(S_n^l).$ Then A has the property $(C_n^l).$
\end{cor}
\par\noindent\textit{Proof:} The claim follows since for any natural numbers $a,b$ and $c$ we have $$\min(a,b-c)\geq\min(a,b)-c.$$

\section{Serre-type properties for modules}
Because the depth of the zero module is not uniformly considered, we first make:
\begin{conv}\label{de0}
The depth of the zero module is infinite, that is $\de_A(0)=\infty.$
\end{conv}
\begin{rem}\label{Cl}
  Regarding the definition of the property $(S_n)$ for modules, we use the definition from \cite[D\'ef. 5.7.2]{EGA} and \cite{BH}. A slightly different one is used in \cite{EG}.
One can learn more about this in \cite{HuWie} and also at the following link\\
https://mathoverflow.net/questions/22228/what-is-serres-condition-s-n-for-sheaves
\end{rem}

\begin{defn}\label{Cnmod}

Let $n$ and $l$ be natural numbers, $A$ a Noetherian ring and $M$ a finitely generated $A$-module. We say that  $M$ has the property $(C_n^l)$ if $$\de(M_\mathfrak{p})\geq\min(\dim(M_\mathfrak{p}),n)-l,\  \text{for any}\ \ \mathfrak{p}\in\spec(A).$$
\end{defn}

\begin{rem}\label{obser}
\par i) One can see immediately that, like in the case of rings, we have that $(C_n^l)\Rightarrow (C_{n-1}^l)$ and that  $(C_n^l)\Rightarrow (C_n^{l+1}).$
\par\noindent  ii) It is also clear that the property $(C_n^l)$ localizes well,  that is if $M$ has $(C_n^l),$ then $M_\mathfrak{p}$ has $(C_n^l)$ for any prime ideal $\mathfrak{p}\in\spec(A).$
\end{rem}

\begin{thm}\label{charcmdm} Let A be a Noetherian ring, M a finitely generated A-module and $l$ a natural number. Then $\cmd(M)\leq l$ if and only if $M$ has $(C_n^l)$ for any $n\in\mathbb{N}.$
\end{thm}
\par\noindent\textit{Proof:} It is similar to the proof of \ref{charcmd}.

\begin{prop}\label{regelm} Let A be a Noetherian ring, M a finitely generated A-module and $l$ a natural number. If M satisfies $(C_n^l)$ and $x$ is a regular element on $M$, then $M/xM$ satisfies $(C_{n-1}^l).$
\end{prop}

\par\noindent\textit{Proof:} Follows at once from the definitions.

\begin{prop}\label{charac} Let A be a Noetherian ring, M a finitely generated A-module and $l$ a natural number. Then  $M$ has $(C_n^l)$ if and only if $\cmd(M_\mathfrak{p})\leq l,$ for any prime ideal $\mathfrak{p}$ such that $\de(M_\mathfrak{p})\leq n-l-1.$
\end{prop}

\par\noindent\textit{Proof:} Let $\mathfrak{p}\in\spec(A).$ If $M_\mathfrak{p}=(0),$ that is $\mathfrak{p}\notin\supp(M),$ then by \ref{de0} we have $\de(M_\mathfrak{p})=\infty.$ This shows that we may assume that $\mathfrak{p}\in\supp(M).$
\par\noindent Assume that $A$ has the property $(C_n^l).$ Let $\mathfrak{p}\in\spec(A)$ with $\de(M_\mathfrak{p})\leq n-l-1.$ Hence $\min(n,\dim(M_\mathfrak{p}))-l\leq\de(M_\mathfrak{p})\leq n-l-1.$
\par\noindent If $\dim(M_\mathfrak{p})\leq n,$ then $\dim(M_\mathfrak{p})-l\leq\de(M_\mathfrak{p}),$ that is $\cmd(M_\mathfrak{p})\leq l.$
\par\noindent If $\dim(M_\mathfrak{p})>n,$ then $n-l\leq\de(M_\mathfrak{p})\leq n-l-1,$ that is obviously a contradiction.
\par\noindent Conversely, let $\mathfrak{p}\in\supp(M).$
\par\noindent If $\de(M_\mathfrak{p})\leq n-l-1,$ then $\cmd(M_\mathfrak{p})=\dim(M_\mathfrak{p})-\de(A_\mathfrak{p})\leq l,$ hence $\dim(M_\mathfrak{p})\leq n-1$. This means that $\min(n,\dim(M_\mathfrak{p}))=\dim(M_\mathfrak{p}),$ that is $\de(M_\mathfrak{p})\leq\min(n,\dim(M_\mathfrak{p}))-l.$
\par\noindent If $\de(M_\mathfrak{p})> n-l-1,$ then $\dim(M_\mathfrak{p})>n-2,$ hence $\de(M_\mathfrak{p})\geq\min(n,\dim(M_\mathfrak{p}))-l.$

\begin{thm}\label{characmod} {\rm (cf. \cite[Th. 2.3]{MT2})}
Let A be a Noetherian ring, M a finitely generated A-module and $l$ a natural number. The following are equivalent:
\par\noindent i) M has the property $(C_n^l)$.
\par\noindent ii) $\ass(M/(x_1,\ldots,x_i)M\subseteq\{\mathfrak{p}\in\supp(M/(x_1,\ldots,x_i)M)\ \vert\ \linebreak \dim(M/(x_1,\ldots,x_i)_\mathfrak{p})\leq l\}$ for any $M$-regular sequence $x_1,\ldots,x_i,$ for any $0\leq i<n-1$.
\par\noindent iii) $\gr(\mathfrak{p},M)\geq\min\{n,\dim(M_\mathfrak{p})\}-l$ for any prime ideal $\mathfrak{p}\in\supp(M).$
\end{thm}

\par\noindent\textit{Proof:} $i)\Rightarrow ii):$ From Prop. \ref{regelm} it follows that $M/(x_1,\ldots,x_i)M$ has the property $(C_{n-i}^l),$ hence $\de((M/(x_1,\ldots,x_i)M_\mathfrak{p})\geq\min\{n-i,\dim(M_\mathfrak{p})\}-l.$ But $\de((M/(x_1,\ldots,x_i)M_\mathfrak{p})=0$ and $n-i>1.$ This forces $\dim(M_\mathfrak{p})\leq l.$
\par\noindent $ii)\Rightarrow iii):$ Assume that there exists a prime ideal $\mathfrak{p}\in\supp(M)$ such that $r:=\gr(\mathfrak{p},M)<\min\{n,\dim(M_\mathfrak{p})\}-l.$ Let $x_1,\ldots,x_r\in \mathfrak{p}$ be a maximal $M$-sequence, where $r<n-1.$ Then $\mathfrak{p}\in\ass(M/(x_1,\ldots,x_r)M$, hence $\dim(M/(x_1,\ldots,x_r)M_\mathfrak{p})\leq l.$ Then $\dim(M_\mathfrak{p})\leq r+l,$ but this contradicts the initial choice of $\mathfrak{p}.$
\par\noindent $iii)\Rightarrow i):$ Let $\mathfrak{p}\in\supp(M).$ Then $\gr(\mathfrak{p},M)\leq\de(M_\mathfrak{p}),$ by \cite[Prop. 1.2.10]{BH}. The conclusion follows.

\begin{prop}\label{deform} {\rm (cf. \cite[Prop. 2.7]{MT2})} Let A be a Noetherian ring, M a finitely generated A-module and $l$ a natural number. Assume that $x$ is an element in the Jacobson radical of $A$ that is a regular element on $M$ and $A$.  If $M/xM$ satisfies $(C_n^l)$ as an $A/xA$-module, then $M$ satisfies $(C_n^l)$ as an $A$-module.
\end{prop}

\par\noindent\textit{Proof:}
Let $\mathfrak{q}\in\spec(A)$ be such that $\de(M_\mathfrak{q})=s\leq n- l-1.$ Note that by \ref{de0} this means that $\mathfrak{q}\in\supp(M).$ If $x\in \mathfrak{q},$ then $\de(M/xM)_\mathfrak{q}=s-1\leq n-l-2.$ Then $\dim(M/xM)_{(\mathfrak{q}/xA)}\leq s-1+1=s,$
hence $\dim(M_\mathfrak{q})\leq s+1=\de (M_\mathfrak{q})+1.$  If $x\notin \mathfrak{q},$ let $\mathfrak{p}\in\smin(\mathfrak{q}+xA).$ Then $(\mathfrak{q}+xA)A_\mathfrak{p}$ is $\mathfrak{p}A_\mathfrak{p}$-primary and $\de(M_\mathfrak{p})\leq\de(M_\mathfrak{q})+1=s+1.$
Then $\de(M/xM)_\mathfrak{q}=s-1,$ hence $\dim(M_{(\mathfrak{p}/xA)})\leq s.$ It follows that $\dim(M_\mathfrak{p})\leq s+1=\de(M_\mathfrak{p})+1.$

\begin{cor}\label{series} {\rm (cf. \cite[Cor. 2.8]{MT2})}
 Let $A$ be a Noetherian ring, $M$ a finitely generated A-module and $l$ a natural number. If the $A$-module $M$ satisfies $(C_n^l)$, then the $A[[X]]-$ module $M[[X]$ satisfies $(C_n^l)$.
\end{cor}

\par \textit{Acnowledgement}. The author is deeply grateful to the referee for the careful reading of the manuscript and for giving useful comments.


\begin{flushright}
 Simion Stoilow Institute of Mathematics \\of the Romanian Academy\\
  PO Box 1-764\\ 014700 Bucharest - Romania\\
  email: cristodor.ionescu@imar.ro\\ cristodor.ionescu@gmail.com
\end{flushright}

  \end{document}